\documentclass{commat}

\usepackage{amscd, tikz-cd}

\newcommand{\z}{\mathbb{Z}}
\newcommand{\Modd}{{\mathcal{M}od}}
\newcommand{\SL}{{\rm SL}}
\newcommand{\GL}{{\rm GL}}
\newcommand{\rr}{{R^\times}}
\newcommand{\mt}{\mapsto}
\newcommand{\se}{\subseteq}
\newcommand{\arr}{\rightarrow}
\newcommand{\tail}{\rightarrowtail}
\newcommand{\larr}{\longrightarrow}
\newcommand{\harr}{\hookrightarrow}
\newcommand{\two}{\twoheadrightarrow}

\DeclareMathOperator{\tor}{Tor}
\DeclareMathOperator{\ext}{Ext}
\DeclareMathOperator{\diag}{diag}
\DeclareMathOperator{\inc}{inc}
\DeclareMathOperator{\corr}{cor}
\DeclareMathOperator{\ress}{res}

\title{%
Some remarks on the homology of nilpotent groups
}

\author{%
    Behrooz Mirzaii and Fatemeh Y. Mokari
    }

\affiliation{
    \address{Behrooz Mirzaii --
  Instituto de Ci\^encias Matem\'aticas e de Computa\c{c}\~ao (ICMC), Universidade de S\~ao Paulo, S\~ao Carlos, Brasil
        }
    \email{%
    bmirzaii@icmc.usp.br
    }
    \address{Fatemeh Yeganeh Mokari --
     Instituto de Ci\^encias Matem\'aticas e de Computa\c{c}\~ao (ICMC), Universidade de S\~ao Paulo, S\~ao Carlos, Brasil
             }
    \email{%
    f.mokari61@gmail.com
    }
    }

\abstract{%
In this article we study the homology of nilpotent groups. In particular,
a certain vanishing result for the homology and cohomology of nilpotent 
groups is proved.}



\VOLUME{31}
\YEAR{2023}
\NUMBER{1}
\firstpage{359}
\DOI{https://doi.org/10.46298/cm.10453}

\begin{paper} 

\section*{Introduction}

The vanishing problem for (co)homology of nilpotent groups asks:
For a nilpotent group $G$ and an $RG$-module $M$, $R$ a commutative ring with $1$, when 
the vanishing of the zero (co)homology of $G$ with coefficients in $M$ 
will result in the vanishing of all (co)homologies of $G$ with coefficients 
in $M$?

Such vanishing results have many interesting applications \cite{robinson1976}, \cite{robinson1987}.
Simple examples show that certain finiteness conditions on $M$ are 
needed \cite{dwyer1975}, \cite{robinson1976}. The most general results in this direction are due to Robinson
\cite[Theorem A and Theorem B]{robinson1976}. He proved that if $M$ is a Noetherian $G$-module 
(resp. an Artinian $G$-module), then we have the vanishing result
for the homology (resp. cohomology) functors.

In this article we study a vanishing result which satisfy different type of finiteness conditions.
As our main result we show that if $H$ is a normal subgroup of 
finite index of a nilpotent group $G$ such that $G/H$ is $l$-torsion and 
if $R$ is a principal ideal domain with $1/l\in R$ and $M$ an $RG$-module, then $M_G=0$ 
implies that $H_n(G, M_H)=0$ for all $n\geq 0$. Similarly $M^G=0$ implies that $H^n(G, M^H)=0$ for all $n\geq 0$.

We generalize these results by removing the conditions $M_G=0$ and $M^G=0$. In fact, we show that for any
$n\geq 0$ the natural maps of pairs
\[
(\inc, \corr):(H, M_G) \arr (G, M_H),\ \ \   (\inc, \ress):(H, M^G) \arr (G, M^H)
\]
 induce
the isomorphisms
\[
H_n(H, M_G) \simeq H_n(G,M_H), \ \ \ \ \ 
H^n(G, M^H) \simeq H^n(H,M^G), 
\]
respectively. More generally, we show that for any integers  $n,r\geq 0$ we have the isomorphisms
\[
H_n(H, H_r(G, M)) \simeq  H_n(G, H_r(H, M)),
\]
\[ 
H^n(G, H^r(H, M))\simeq  H^n(H, H^r(G, M)).
\]

\section{A vanishing result}

\begin{theorem}\label{vanishing2}
Let $G$ be a nilpotent group and  $H$ a normal subgroup of $G$ such that 
$G/H$ is finite and $l$-torsion. Let $R$ be a principal ideal domain with $l\in R^\times$ and $M$ an $RG$-module.
\par {\rm (i)} If $M_G=0$, then for any $n\geq 0$, $H_n(G, M_H)=0$.
\par {\rm (ii)} If $M^G=0$, then for any $n\geq 0$, $H^n(G, M^H)=0$.

In particular, if $H$ acts trivially on $M$ and $M_G=0$ (resp. $M^G=0$),
then for any $n\geq 0$,  $H_n(G,M)=0$ (resp. $H^n(G,M)=0$).
\end{theorem}
\begin{proof}
We prove the theorem in few steps:
\par {\bf Step 1.} {\sf If $G$ is finite and $l$-torsion, then the coinvariant and invariant
functors 
\[
(-)_G,\ (-)^G : \Modd_{RG} \arr \Modd_R
\]
are exact:}
Set $F:=()_G$ and $F':=()^G$. First let $l=|G|$.
Since $F$ is right exact and $F'$ is left exact, to prove the claims
it is sufficient to prove that the natural map $\alpha_G: M^G \arr M_G$ 
is an isomorphism. If 
\[
\overline{N}: M_G \arr M^G, \ \ \ \overline{m} \mt Nm,
\]
where $N:=\sum_{g \in G}g \in R G$, then clearly $\overline{N}\circ\alpha_G$ 
and $\alpha_G \circ\overline{N}$ coincide with multiplication by $|G|$.
Thus $\alpha_G$ is an isomorphism. The proof of the general case
is by induction on the size of $G$ and we 
may assume that $G\neq 1$. Since $G$ is nilpotent, $Z(G)\neq 1$.
Let $H$ be a nontrivial cyclic subgroup of $Z(G)$.
The map $\alpha_G$ coincides with the following composition of maps
\[
M^G\overset{\simeq}{\larr} {(M^H)}^{G/H} 
\overset{\alpha_H}{\larr} {(M_H)}^{G/H}
\overset{\alpha_{G/H}}{\larr} {(M_H)}_{G/H}
\overset{\simeq}{\larr} M_G.
\]
Now the claims follow from the above argument and the induction process.

{\bf Step 2.} {\sf If $G$ is $l$-torsion, then
$H_n(G, M)=H^n(G,M)=0$:} 
This is a known fact. But here we give a direct proof of it.  
First let $G$ be finite and $l=|G|$.
If $P_\bullet \arr M$  and $M \arr I^\bullet$ are projective  and
injective resolutions of the $RG$-module $M$, respectively, then the claim 
follows from Step~1 and the following isomorphisms
\[
H_n(G, M)\simeq H_n({(P_\bullet)}_G), \ \ \ 
H^n(G, M) \simeq H_n({(I^\bullet)}^G)
\]
(see \cite[Chap. III.6, 1.4 and III.6, Exercise 1]{brown1994}).
In general, since $G$ is nilpotent and torsion, 
any finitely generated subgroup of $G$ is finite. Hence we can write $G$ 
as direct limit of its finite subgroups, e.g. 
$G=\underset{\larr}{\lim} \ G_i$, where $G_i$'s are finite. For the 
homology functor we have
\[
H_n(G, M)\simeq \underset{\larr}{\lim}\ H_n(G_i, M)=0
\]
(see \cite[Chap. V.5, Exercise 3]{brown1994}). For the 
cohomology functor we have the spectral sequence
\[
E_2^{p,q}=\underset{\longleftarrow}{\lim}^p\ H^q(G_i, M) 
\Rightarrow H^{p+q}(G, M),
\]
where $\underset{\longleftarrow}{\lim}^p$ is the $p$-th derived functor of 
$\underset{\longleftarrow}{\lim}$ \cite[p. 297]{robinson1987}.
Now the claim follows from the finite case.

{\bf Step 3.}
{\sf If $G$ is finite and $l$-torsion, then for any $R$-module 
$N$ with the trivial action of $G$ and any $n \ge 0$, we have the isomorphisms
\[
\tor_n^{R}(N,M)_G \simeq \tor_n^{R}(N,M_G),\ \ \
\ext_R^n(N,M^G) \simeq \ext_R^n(N,M)^G:
\]}
We start with the functor $\tor$. Since
\[\tor_0^{R}(N,M)_G \simeq
(N \otimes_R M) \otimes_G \z
\simeq N \otimes_R M_G  \simeq  \tor_0^{R}(N,M_G),
\]
the claim is true for $n=0$. Let 
$$0 \arr N_{n-1} \arr F \arr N \arr 0$$ 
be a
short exact sequence of $R$-modules such that
$F$ is free. If $ n \ge 2$, from the long exact sequence, we get
the isomorphism $\tor_n^{R}(N,M) \simeq \tor_{n-1}^{R}(N_{n-1},M)$.
If we continue this process, we will find an $R$-module $N_1$ such that
\[
\tor_n^{R}(N,M) \simeq \tor_{1}^{R}(N_{1},M).
\]
So it is sufficient to proof the claim for $n=1$. From the exact sequence
$$0 \arr N_1 \arr F \arr N \arr 0$$ 
and Step 1 we obtain the following commutative 
diagram with exact rows
\[
\begin{CD}
0 \arr & \tor_1^{R}(N ,M)_G  & \larr & (N_1\otimes_R M)_G
& \larr & (F \otimes_R M)_G & \larr
& (N \otimes_R M)_G   & \arr 0 \\
&  @VVV  @VVV @VVV  @VVV & \\
0 \arr & \tor_1^{R}(N ,M_G)  & \larr & N_1\otimes_R M_G
& \larr & F \otimes_R M_G & \larr
& N \otimes_R M_G   & \arr 0.
\end{CD}
\]
Note that the three vertical maps on the right are isomorphisms. Now the claim follows from an easy diagram chase. 

The proof of the claim for the functor  $\ext$ is similar. In fact here we should use an 
exact sequence $0 \arr N \arr I \arr N^1 \arr 0$,  where $I$ is an 
injective $R$-module.

{\bf Step 4.} {\sf 
$H_n(H, M)_{G/H}\simeq H_n(G,M)$ and $H^n(G,M)\simeq H^n(H, M)^{G/H}$:}
We prove the first isomorphism.  The second one can be proved in a similar way.
From the extension $H \tail G \two G/H$ we obtain the Lyndon-Hochschild-Serre homology spectral sequence
\[
E_{p,q}^2=H_p(G/H, H_q(H, M))\Rightarrow H_{p+q}(G, M).
\]
By Step 2, we have $E_{p,q}^2=0$. Now by an easy analysis of the spectral sequence we obtain the isomorphism
$H_0(G/H, H_n(H, M))\simeq H_n(G, M)$.

{\bf Step 5.} {\sf The proof of the theorem:}
We prove (i). The proof of (ii) is similar. By Step~4,
we have the isomorphism 
\[
H_n(G,M_H)\simeq H_n(H,M_H)_{G/H}.
\]
Since the action of $H$ on $M_H$ is trivial, by the Universal Coefficient Theorem 
we have the exact sequence
\[
0 \arr H_n(H,R)\otimes_R M_H \arr H_n(H,M_H) \arr 
\tor_1^R(H_{n-1}(H,R), M_H) \arr 0.
\] 
By applying the functor $(.)_{G/H}$ we obtain the exact sequence
\[
0 \arr\big(H_n(H,R)\otimes_R M_H\big)_{G/H}\arr H_n(H,M_H)_{G/H} 
\arr \tor_1^R(H_{n-1}(H,R), M_H)_{G/H} \arr 0.
\]
Now the claim follows from Step 3 and the assumption $M_G=0$.
\end{proof}

\section{The corestriction map for the homology of nilpotent groups}

We say that a group $G$ acts nilpotently on a $G$-module $M$, if $M$ has a finite filtration 
of $G$-submodules 
\[
0= M_0 \se \cdots \se M_{k-1} \se M_k=M,
\]
such that the action of $G$ on each quotient $M_i/M_{i-1}$ is trivial. 

The following theorem will be needed in the next section.

\begin{theorem}\label{nil-corestriction}
Let $G$ be a nilpotent group and $H$ a normal subgroup of $G$ such that $G/H$ is $l$-torsion. Let $R$ be a 
commutative ring such that $l \in R^\times$. If $M$ is an $R$-module with a nilpotent action of $G$, then for any 
$n\geq 0$ the natural maps
\[
\corr_H^G:H_n(H,M)\arr H_n(G, M), \ \ \ress_H^G:H^n(G,M)\arr H^n(H,M)
\]
are isomorphisms. In particular,  $G/H$ acts trivially on $H_n(H, M)$ and $H^n(H, M)$.
\end{theorem}
\begin{proof}
The proof is by induction on the nilpotent class $c$ of $G$. We prove the claim for $\corr_H^G$. 
The  claim for $\ress_H^G$ can be proved in a similar way. Consider the lower central series of $G$:
\[
1=\gamma_{c+1}(G)\subset \gamma_{c}(G)\subset \cdots \subset 
\gamma_{2}(G)\subset \gamma_{1}(G)=G.
\]

Let $c=1$. Then $G'=\gamma_{2}(G)=1$.  So $G$ is abelian. Let
\[
0=M_0 \se M_1 \se \cdots \se M_k=M,
\] 
be a filtration of $M$ such that $G$ acts trivially on each quotient $M_i/M_{i-1}$. We prove this case by induction on $k$.
If $k=1$, then the action of $G$ on $M=M_1$ is  trivial. This implies that the action of $G/H$ on $H_n(G,M)$ is trivial  and 
therefore
\[
H_n(H,M)= H_n(H,M)_{G/H}\simeq H_n(G,M)
\]
(see Step 4, in the proof of Theorem \ref{vanishing2}).
Now let $k>1$ and set $M_1':=M/M_1$. From the  short exact sequence of $G$-modules 
$0 \arr M_1 \arr M \arr M_1' \arr 0$, we obtain the following commutative diagram with exact rows
\[
\begin{CD}
H_{n+1}(H, M_1') &\arr & H_n(H, M_1) &\arr &H_n(H, M) &\arr& H_n(H, M_1') &\arr& H_{n-1}(H, M_1)\\
  @V{f_1} VV                      @V{g_1}VV               @V{\corr_H^G}VV              @V{f_2}VV                    @V{g_2}VV\\
H_{n+1}(G, M_1') &\arr&  H_n(G, M_1) &\arr & H_n(G, M) &\arr&  H_n(G, M_1') &\arr& H_{n-1}(G, M_1),
\end{CD}
\]
where $f_1$, $f_2$, $g_1$ and $g_2$ are the natural corestriction maps. Since $G$ acts trivially on $M_1$ and $M_1'$ has a filtration 
of length $k-1$, by induction $f_1$, $f_2$, $g_1$ and $g_2$ are isomorphisms. Now an easy diagram chase shows that $\corr_H^G$ is 
an isomorphism. This proves the theorem for $c=1$.

Now assume that the claim is true for nilpotent groups of class $d$,  $1\leq d \leq c-1$. From the commutative diagram of extensions,
\[
\begin{tikzcd}
\gamma_{c}(G) \cap H \ar[r, tail]\ar[d]&  H \ar[r, two heads]\ar[d] &H/(\gamma_{c}(G) \cap H ) \ar[d] \\ 
\gamma_{c}(G)   \ar[r, tail]     &  G  \ar[r, two heads]&  G/\gamma_{c}(G),    
\end{tikzcd}
\]
we have the following morphism of Lyndon-Hochschild-Serre spectral sequences
\[
\begin{tikzcd}
E'^2_{p,q}=H_p(H/(\gamma_{c}(G)\cap H),H_q(\gamma_{c}(G)\cap H,M))\ar[r, Rightarrow] \ar[d]& H_{p+q}(H, M)\ar[d]\\ 
E^2_{p,q}=H_p(G/\gamma_{c}(G) ,H_q(\gamma_{c}(G) , M)) \ar[r, Rightarrow] & H_{p+q}(G, M).
\end{tikzcd}
\] 
First note that 
$\gamma_{c}(G)/(\gamma_{c}(G) \cap H)\simeq  \gamma_c(G)H/H \se G/H$. So the groups 
$\gamma_{c}(G)/(\gamma_{c}(G) \cap H)$ and $(G/\gamma_{c}(G))/(H/(\gamma_{c}(G)\cap H))$ are $l$-torsion.

Since $\gamma_c(G)$ is abelian, by the first step of the induction, we have
\[
H_q(\gamma_{c}(G) \cap H, M)\simeq H_q(\gamma_{c}(G), M). 
\]
Observe that $G/\gamma_{c}(G)$ is of nilpotent class $c-1$.
Since $\gamma_{c}(G)\se Z(G)$, the conjugate action of $G/\gamma_{c}(G)$ on $\gamma_{c}(G)$ is trivial. 

We show that the natural action of $G/\gamma_{c}(G)$ on $H_q(\gamma_{c}(G), M)$  is nilpotent. 
This can be done by induction on the length of the filtration of $M$. Assume that $M$ has a filtration of length $k$ as above. 
If $k=1$, then $M=M_1$. Thus $G$ acts trivially on $M$ and so the action of $G/\gamma_{c}(G)$ on $H_q(\gamma_{c}(G), M)$  
is trivial. Now let $k\geq 2$. From the short exact sequence 
$$0\arr M_1 \arr M \arr M/M_1 \arr 0,$$
we obtain the long exact sequence
\[
\cdots \arr H_q(\gamma_{c}(G), M_1)\arr H_q(\gamma_{c}(G), M)\arr H_q(\gamma_{c}(G), M/M_1)\arr \cdots.
\]
By induction, the actions of $G/\gamma_{c}(G)$ on the  modules $H_q(\gamma_{c}(G), M_1)$ and $H_q(\gamma_{c}(G), M/M_1)$ 
are nilpotent. It follows from the above exact sequence that $G/\gamma_{c}(G)$ acts nilpotently on $H_q(\gamma_{c}(G), M)$.

Since $H/(\gamma_{c}(G) \cap H)\harr G/\gamma_{c}(G)$ and $H_q(\gamma_{c}(G), M)\simeq H_q(\gamma_{c}(G)\cap H, M)$,
by  induction on the nilpotent class of $G/\gamma_{c}(G)$, we have
\[
H_p(H/(\gamma_{c}(G) \cap H ), H_q(\gamma_{c}(G)\cap H, M))\simeq H_p(G/\gamma_{c}(G), H_q(\gamma_{c}(G), M)).
\]
Therefore $E'^2_{p,q}\simeq E^2_{p,q}$. Now by convergence of the spectral sequences, for any $n\geq 0$, we obtain the isomorphism
\[
H_n(H, M)\simeq H_n(G, M).
\]
Finally we know that $H_n(H, M)_{G/H} \simeq H_n(G, M)$ (see Step 4, in the proof of Theorem \ref{vanishing2})
Thus the map $H_n(H,M)\arr H_n(H, M)_{G/H}$ is an isomorphism. This shows that $G/H$ acts trivially on $H_n(H,M)$.
\end{proof}

\begin{example}
Easy examples show that in Theorem \ref{nil-corestriction} the condition 
$l\in R^\times$ can not be removed. For example, if $H$ is a proper finite subgroups of an abelian group $G$, 
then $H_1(H, \z)\simeq H \subset G\simeq H_1(G,\z)$, which clearly is not an isomorphism.
\end{example}

The first part of Theorem \ref{nil-corestriction} can be extended to all subgroups 
of finite index, as follows.

\begin{corollary}\label{nil-transfer}
Let $G$ be a nilpotent group and $H$ a subgroup of finite index. Let $R$ be 
a commutative ring such that $[G:H]! \in R^\times$. If $M$ is an $R$-module with a
nilpotent action of $G$, then, for any $n\geq 0$, the natural maps
\[
\corr_H^G:H_n(H,M)\arr H_n(G, M), \ \ \ress_H^G:H^n(G,M)\arr H^n(H,M)
\]
are isomorphisms. 
\end{corollary}

\begin{proof}
It is well-known that $H$ has a subgroup $L$ such that $L$ is normal in $G$ and 
\[
[G:L] \leq [G:H]!.
\]
By Theorem \ref{nil-corestriction}, $\corr_L^G:H_n(L,M)\arr H_n(G, M)$ and $\corr_L^H:H_n(L,M)\arr H_n(H, M)$
are isomorphisms. Therefore 
\[
\corr_H^G:H_n(H,M)\arr H_n(G, M)
\]
is an isomorphism. The cohomology case can be treated similarly.
\end{proof}

\section{Homology of nilpotent groups with corestriction coefficients }

The following theorem generalizes Theorem \ref{vanishing2}.

\begin{theorem}\label{a-b}
Let $G$ be a nilpotent group and  $H$ a normal subgroup of $G$ such that $G/H$ 
is finite and $l$-torsion. Let $R$ be a principal ideal domain with $l \in R^\times$ and let $M$ be an 
$RG$-module. Then for any $n \ge 0$,
\[
H_n(H, M_G) \simeq  H_n(G, M_H), \ \ \ \ H^n(G, M^H)\simeq  H^n(H, M^G),
\]
which are induced by the pairs
\[
    (\inc, \corr):(H, M_G) \arr (G, M_H)
    \qquad \text{and} \qquad
    (\inc, \ress):(H, M^G) \arr (G, M^H),
\]
respectively. More generally,
for any integers $n,r \geq 0$, the above maps of pairs induce the isomorphisms
\[
H_n(H, H_r(G, M)) \simeq  H_n(G, H_r(H, M)),
\]
\[ 
H^n(G, H^r(H, M))\simeq  H^n(H, H^r(G, M)).
\]
\end{theorem}
\begin{proof}
The Universal Coefficient Theorem 
and Step 1 of the proof of Theorem \ref{vanishing2} gives us the
following commutative diagram with exact rows
\[
\begin{CD}
0 \arr & (H_n(H, R) \otimes_R M_H)_{G/H}  & \arr & H_n(H, M_H)_{G/H} & \arr &
\tor_1^R(H_{n-1}(H, R), M_H)_{G/H} & \arr 0 \\
        &  @VVV  @VVV @VVV & \\
0 \arr & H_n(H, R)\otimes_R M_G  & \arr & H_n(H, M_G) & \arr &
\tor_1^R(H_{n-1}(H, R), M_G) & \arr 0.
\end{CD}
\]
By Theorem \ref{nil-corestriction} the action of $G/H$ on $H_n(H,R)$ 
is trivial. Thus by Step 3 of the proof of Theorem \ref{vanishing2}, 
the left and the right column maps of this diagram are isomorphisms. Hence 
$H_n(H, M_H)_{G/H} \simeq  H_n(H, M_G)$ and therefore
\[
H_n(G, M_H)\simeq H_n(H, M_H)_{G/H} \simeq H_n(H, M_G).
\]
The other isomorphism can be proved in a 
similar way. In fact, one should prove that 
$H^n(H, M^G)\simeq H^n(H, M^H)^{G/H}$ first, and then use this result to prove the isomorphism $H^n(H, M^G) \simeq H^n(G, M^H)$.

The proof of the general case is similar. Just we should replace $M_H$ and $M_G$, with
$H_r(H,M)$ and $H_r(G,M)$, respectively.
\end{proof}

\begin{corollary}\label{a-bb}
Let $G$ be a nilpotent group, $H$ be a subgroup of $G$ such that $G/H$ is
$l$-torsion, $R=\z[1/l]$, and $M$ be an $RG$-module. Then for any $n \ge 0$,
$H_n(G, M_H) \simeq  H_n(H, M_G)$. In particular, if the action of $H$ on $M$ is
trivial, then
\[
H_n(G, M) \simeq  H_n(H, M_G).
\]
\end{corollary}
\begin{proof}
First notice that the group $G/H$ can be written as direct limit of its finite subgroups, e.g. $G/H=\underset{\larr}{\lim}\ G_i/H$. (Note that finitely generated torsion subgroups of nilpotent groups are finite.) Hence by Theorem \ref{a-b},
\[
H_n(G, M_H)\simeq \underset{\larr}{\lim}\ H_n(G_i, M_{H})
\simeq \underset{\larr}{\lim}\ H_n(H, M_{G_i})
\]
\[
\hspace{1.5 cm}
\simeq H_n(H, \underset{\larr}{\lim}\ M_{G_i})\simeq H_n(H,M_G).
\]
(see \cite[Exercise 3, Chapter V.5]{brown1994}).
\end{proof}

\begin{example}
As an application of Theorem \ref{a-b}, we study the homology 
of special linear groups. Let $R$ be a commutative ring. The conjugate 
action of $\rr$ on $\SL_n(R)$, given by
\[
a.A:=\diag(a, I_{n-1}).A.\diag(a^{-1}, I_{n-1}),
\]
induces a natural action of $\rr$ on $H_q(\SL_n(R), \z)$. Since
\[
\hspace{-4.4cm}
a^n.A=\diag(a^n, I_{n-1}).A.\diag(a^{-n}, I_{n-1})
\]
\[
=\diag(a^{n-1}, a^{-1}I_{n-1}).aI_n.A.a^{-1}I_n.\diag(a^{-(n-1)}, aI_{n-1})
\]
\[
\hspace{-1.6cm}
=\diag(a^{n-1}, a^{-1}I_{n-1}).A.\diag(a^{-(n-1)}, aI_{n-1}),
\]
the action of $\rr^n$ on
$H_q(\SL_n(R), \z)$ is trivial \cite[Chap. II, Proposition 6.2]{brown1994}.
Since $\rr/\rr^n$ is a $n$-torsion group, by Corollary \ref{a-bb}
\[
H_p(\rr, H_q(\SL_n(R), \z[1/n])) \overset{\simeq}{\larr} H_p(\rr, H_q(\SL_n(R), \z[1/n])_\rr).
\]

We say that a commutative ring $R$ is a {\it ring with many units}
if for any $n \ge 2$ and for any finite number of surjective linear
forms $f_i: R^n \arr R$, there exists a $v \in  R^n$ such that, for
all $i$, $f_i(v) \in \rr$. Important examples of rings with many units
are semi-local rings with infinite residue fields. For more about these
rings please see \cite[Section 1]{guin1989} and \cite[Section 2]{mirzaii2008}.
Now one can show that if $q \le n$ are nonnegative integers, then
\[
H_q(\SL_n(R), \z[1/n])_\rr \simeq H_q(\SL(R),\z[1/n]), 
\]
provided that $R$ is a ring with many units \cite[Section 3]{mirzaii2008}. Combining these results, 
we obtain the isomorphism
\[
H_p(\rr, H_q(\SL_n(R),\z[1/n]))\overset{\simeq}{\larr} 
H_p(\rr, H_q(\SL(R),\z[1/n])).
\]
\end{example}

\subsection*{Acknowledgments} 
The second author  acknowledges that during the work on this paper she was supported by a post-doc fellowship 
of FAPESP (Funda\c{c}\~ao de Amparo \`a Pesquisa do Estado de S\~ao Paulo) with grant number 2016/13937-9.


\EditInfo{August 04, 2020}{January 13, 2022}{Ivan Kaygorodov}

\end{paper}